\documentstyle[12pt]{article}

\def\bea{\begin{eqnarray}}
\def\eea{\end{eqnarray}}

\def\beq{\begin{equation}}
\def\eeq{\end{equation}}
\def\ba{\beq\new\begin{array}{c}}
\def\ea{\end{array}\eeq}
\def\be{\ba}
\def\ee{\ea}

\makeatletter
\newdimen\normalarrayskip 
\newdimen\minarrayskip 
\normalarrayskip\baselineskip
\minarrayskip\jot
\newif\ifold \oldtrue \def\new{\oldfalse}
\def\arraymode{\ifold\relax\else\displaystyle\fi} 
\def\eqnumphantom{\phantom{(\theequation)}} 
\def\@arrayskip{\ifold\baselineskip\z@\lineskip\z@
\else
\baselineskip\minarrayskip\lineskip2\minarrayskip\fi}
\def\@arrayclassz{\ifcase \@lastchclass \@acolampacol \or
\@ampacol \or \or \or \@addamp \or
\@acolampacol \or \@firstampfalse \@acol \fi
\edef\@preamble{\@preamble
\ifcase \@chnum
\hfil$\relax\arraymode\@sharp$\hfil
\or $\relax\arraymode\@sharp$\hfil
\or \hfil$\relax\arraymode\@sharp$\fi}}
\def\@array[#1]#2{\setbox\@arstrutbox=\hbox{\vrule
height\arraystretch \ht\strutbox
depth\arraystretch \dp\strutbox
width\z@}\@mkpream{#2}\edef\@preamble{\halign
\noexpand\@halignto
\bgroup \tabskip\z@ \@arstrut \@preamble \tabskip\z@ \cr}%
\let\@startpbox\@@startpbox \let\@endpbox\@@endpbox
\if #1t\vtop \else \if#1b\vbox \else \vcenter \fi\fi
\bgroup \let\par\relax
\let\@sharp##\let\protect\relax
\@arrayskip\@preamble}
\def\eqnarray{\stepcounter{equation}%
\let\@currentlabel=\theequation
\global\@eqnswtrue
\global\@eqcnt\z@
\tabskip\@centering
\let\\=\@eqncr
$$%
\halign to \displaywidth\bgroup
\eqnumphantom\@eqnsel\hskip\@centering
$\displaystyle \tabskip\z@ {##}$%
\global\@eqcnt\@ne \hskip 2\arraycolsep
$\displaystyle\arraymode{##}$\hfil
\global\@eqcnt\tw@ \hskip 2\arraycolsep
$\displaystyle\tabskip\z@{##}$\hfil
\tabskip\@centering
&{##}\tabskip\z@\cr}
\begingroup\ifx\undefined\newsymbol \else\def\input#1 {\endgroup}\fi

\def\nl#1#2{\mathop{#2}\limits_{#1}}

\newcommand{\pl}{\partial}

\newcommand{\ol}{\overline}

\begin{document}

\setcounter{footnote}{1}
\def\thefootnote{\fnsymbol{footnote}}
\begin{center}
\vspace{0.3in}
{\Large\bf On projective group properties}\\
\smallskip
{\Large\bf of the $6D$ $H$-space of the type $[33]$}
\end{center}
\smallskip
\centerline{{\large
Zolfira Zakirova}\footnote{Kazan State University,
Kazan, Tatarstan, Russia;
e-mail: zolya\_zakirova@mail.ru
}}

\bigskip

\abstract{\small  In this note we find the metric of 6-dimensional
$h$-space of the $[33]$ type and then determine an important
projective group characteristic of this $h$-space.
}

\begin{center}
\rule{5cm}{1pt}
\end{center}

\bigskip
\setcounter{footnote}{0}
\section{Introduction}

The problem of defining $2D$ Riemannian manifolds which admit
projective motions, i. e. continuous transformation groups preserving
geodesics dates back to  S.Lie who considered it in \cite{sof}.
In more recent times, that was
A.V.Aminova who has got a complete solution to this problem \cite{am1}.
For the Riemannian manifolds of dimension greater than 2
the same problem has been solved
by G.Fubini \cite{fub} and  A.S.Solodovnikov \cite{sol}.

In the later paper, \cite{am2} A.V.Aminova has classified
all the Lorentzian manifolds
of dimension $\geq 3$ that admit nonhomothetic projective or affine
infinitesimal transformations. In each case, there were determined the
corresponding maximal projective and affine Lie algebras.

In this paper we continue to investigate this problem and
study the $6$-dimensional pseudo-Riemannian space $V^6(g_{ij})$
with signature
$[++----]$, and, particularly, the $h$-space of the type $[33]$.
The general method of determining pseudo-Riemannian manifolds that
admit some nonhomothetic projective group $G_r$ has been developed
by A.V.Aminova \cite{am2} and has got a detailed account in the paper
by the author \cite{zak} for the $h$-space of the type $[51]$.

\section{The metric of the type $[33]$ $h$-space}

In order to find the $h$-space admitting a nonhomothetic
infinitesimal projective transformation, one needs to integrate the
Eisenhart equation
\be\label{0.1}
h_{ij,k}=2g_{ij} \varphi_{,k}+g_{ik} \varphi_{,j}+
g_{jk} \varphi_{,i},
\ee
which, in the skew-normal frames, is of the form \cite{am2}
\be\label{1}
X_r \overline a_{pq}+\sum_{h=1}^n e_h( \overline a_{hq} \gamma _{\tilde hpr}+
\overline a_{ph} \gamma_{\tilde  hqr})=\overline g_{pr}X_q \varphi+\overline
g_{qr}X_p \varphi
\quad
(p, q, r=1,\ldots,n),
\ee
where $X_r \varphi \equiv {\nl{r}\xi}^i \frac{\partial \varphi}
{\partial x^i}$,   $\gamma_{pqr}=-\gamma_{qpr}=
{\nl{p}{\xi}\hspace{-2.5mm}\phantom{a}_{i,j}}{\nl{q}\xi}^i{\nl{r}\xi}^j$,
${\nl{i}\xi}^j$ are the components of skew-normal frames,
$\overline g_{pr}$ and
$\overline a_{pq}$ are the canonical forms of the tensors
$g_{pr}$, $a_{pq}$, respectively. For the type $[33]$  $h$-space
these latter have the following form \cite{pet}
\be\label{2}
\ol{g}_{ij} dx^i dx^j=e_3 (2 dx^1 dx^3+{dx^2}^2) +
e_6 (2 dx^4 dx^6+{dx^5}^2),
\ee
$$
\ol{a}_{ij} dx^i dx^j=e_3 \lambda_3(2 dx^1 dx^3+{dx^2}^2)+
2 e_3 dx^2 dx^3+
e_6 \lambda_6(2 dx^4dx^5+{dx^5}^2)+
2 e_6 dx^5 dx^6,
$$
where $e_3, e_6 =\pm 1$ and $\lambda_3, \lambda_6$ are the roots of the
characteristic equation ${\det}(h_{ij}-{\lambda}g_{ij})=0$.
Here  $h_{ij}={\nl{\xi}L} g_{ij}$, the Lie derivative of the metric
$g_{ij}$ in the direction of vector ${\xi}^i$ gives the projective
motion  $x^{i'}=x^i+{\xi}^i \delta t$.

Substituting in (\ref{1}) the canonical forms $\overline g_{pr}$ and
$\overline a_{pq}$ from (\ref{2}) and taking into account that
for the type $[33]$ $h$-space $e_1=e_2=e_3$, $e_4=e_5=e_6$,
$\tilde1=3$,  $\tilde2=2$,  $\tilde3=1$,
$\tilde4=6$,  $\tilde5=5$,  $\tilde6=4$, one gets the following
system of equations
\be\label{3}
X_r \lambda_3=0 \; (r \ne 3),
\quad
X_r \lambda_6=0 \; (r \ne 6),
\\
X_3(\lambda_3-2/3 \varphi)=X_6 (\lambda_6 -2/3 \varphi)=0,
\\
\gamma_{213}=1/2 e_3 X_3 \lambda_3,
\quad
\gamma_{312}=\gamma_{321}=-e_3 X_3 \varphi,
\\
\gamma_{546}=1/2 e_6 X_6 \lambda_6,
\quad
\gamma_{645}=\gamma_{654}=-e_6 X_6 \varphi,
\\
\gamma_{163}=\gamma_{262}= \gamma_{361}= \frac{e_3 X_6{\varphi}}
{\lambda_3-\lambda_6},
\quad
\gamma_{263}=\gamma_{362}= -\frac{e_3 X_6 \varphi}{(\lambda_3-\lambda_6)^2},
\\
\gamma_{346}=\gamma_{355}= \gamma_{364}= \frac{e_6 X_3{\varphi}}
{\lambda_3-\lambda_6},
\quad
\gamma_{536}=\gamma_{635}= -\frac{e_6 X_3 \varphi}{(\lambda_3-\lambda_6)^2},
\\
\gamma_{363}=\frac{e_3 X_6 \varphi}{(\lambda_3-\lambda_6)^3},
\quad
\gamma_{366}=\frac{e_6  X_3 \varphi}{(\lambda_3-\lambda_6)^3}.
\ee
All others $\gamma_{pqr}$ are equal to zero.

The system of linear partial differential equations
$$
X_q \theta={\nl{q}\xi}^i\pl_i \theta=0,
\quad
(q=1,\ldots,m, i=1,\ldots,6),
$$
where ${\nl{q}\xi}^i$ are the components of skew-normal frames, is
completely integrable if and only if the following conditions are fulfilled
(see \cite{ezen1}, p. 143)
\be\label{4}
(X_q, X_r) \theta=X_q X_r \theta-X_r X_q \theta=
\sum_{p=1}^6 e_p(\gamma_{pqr}-\gamma_{prq}) X_{\tilde p}.
\ee
Using formulas (\ref{3}) and (\ref{4}), one can write
the commutators of the
operators $X_i$ $(i=1,\ldots,6)$ in the $h$-space under consideration,
\be\label{5}
(X_1,X_2)=(X_1,X_4)=(X_1,X_5)=(X_2,X_5)=(X_4,X_5)=0,
\\
(X_1,X_3)=e_2(\gamma_{213}-\gamma_{231})X_2,
\quad
(X_1,X_6)=-e_3\gamma_{361}X_1,
\\
(X_2,X_3)=e_1(\gamma_{123}-\gamma_{132})X_3,
\quad
(X_2,X_6)=-e_3 \gamma_{362} X_1-e_2\gamma_{262}X_2,
\\
(X_3,X_4)=e_6 \gamma_{634}X_4,
\quad
(X_3,X_5)=e_6\gamma_{635}X_4+e_5 \gamma_{535} X_5,
\\
(X_3,X_6)=-e_3 \gamma_{363}X_1-e_2 \gamma_{263} X_2-
e_1 \gamma_{163}X_3+e_6 \gamma_{636}X_4+e_5\gamma_{536} X_5+e_4
\gamma_{436}X_6,
\\
(X_4,X_6)=e_5(\gamma_{546}- \gamma_{564}) X_5,
\quad
(X_5,X_6)=e_4 (\gamma_{456}- \gamma_{465}) X_6.
\ee
From these relations it follows that
the systems $X_p \theta=0$ $(p \ne 3)$ and
$X_r \theta=0$ $(r \ne 6)$ are completely integrable, their
solutions being $\theta^3$  and $\theta^6$ respectively. The systems
$X_1 \theta=X_{\alpha} \theta=0$ $(\alpha=4, 5, 6)$,
$X_4 \theta=X_{\beta} \theta=0$ $(\beta=1, 2, 3)$,
$X_{\alpha} \theta=0$ and $X_{\beta} \theta=0$  are also integrable.
The first system have solutions $\theta^2$, $\theta^3$, the
second system -- $\theta^5$ and $\theta^6$, the third system --
$\theta^1$, $\theta^2$, $\theta^3$ and fourth system have solutions
$\theta^4$, $\theta^5$ and  $\theta^6$. After coordinate transformation
$x^{i'}=\theta^i(x)$, one gets
\be\label{6}
{\nl{s}\xi}^i=P_s(x){\delta_s}^i,
\quad
{\nl{2}\xi}^3={\nl{2}\xi}^{\alpha}={\nl{3}\xi}^{\alpha}=
{\nl{5}\xi}^{\beta}={\nl{5}\xi}^6={\nl{6}\xi}^{\beta}=0,
\quad
(s=1, 4).
\ee
Now integrating the system of equations (\ref{5}) with using
(\ref{3}) and (\ref{6}), one can find the components
${\nl{i}\xi}^j$ of skew-normal frames.
Using these results and the formulas \cite{am2}
\be\label{7}
g^{ij}=\sum_{h=1}^n e_h {\nl{h}{\xi}}^i {\nl{\tilde h}{\xi}}^j,
\quad
{\nl{h}{\xi}\hspace{-2.5mm}\phantom{a}_{i}}=g_{ij} {\nl{h}\xi}^j,
\quad
a_{ij}=\sum_{h, l=1}^n {e_h e_l} {\ol a}_{hl}
{\nl{\tilde h}{\xi}\hspace{-2.5mm}\phantom{a}_{i}} {\nl{\tilde l}{\xi}\hspace{-2.5mm}\phantom{a}_{j}}
\ee
we can calculate the components of tensors $g_{ij}$ and $a_{ij}$.
Then we come to the following

\bigskip

\noindent
{\bf Theorem 1.} {\it If a symmetric tensor $h_{ij}$ of
the characteristics ${\rm [33]}$ and a function $\varphi$
satisfy in $V^6(g_{ij})$
equation {\rm (\ref{0.1})}, then there exists a holonomic coordinate system
so that the function $\varphi$ and the tensors $g_{ij}$,
$h_{ij}$ are defined by the formulas
\be\label{8}
g_{ij}dx^idx^j=e_3 (f_6-f_3)^3 \lbrace (dx^2)^2 +4A  dx^1dx^3+
2(\epsilon x^1-2 A \Sigma_1)dx^2 dx^3+
\ee
$$
((\epsilon x^1)^2-4\epsilon x^1 A \Sigma_1+4A^2 \Sigma_2) (dx^3)^2 \rbrace+
e_6 (f_3-f_6)^3
\lbrace (dx^5)^2+4 \tilde{A}  dx^4dx^6+
2(\tilde{\epsilon} x^4+
$$
$$
2 \tilde{A} {\Sigma}_1)dx^5dx^6+
((\tilde{\epsilon} x^4)^2+4\tilde{\epsilon} x^4 \tilde{A} {\Sigma}_1+4
\tilde{A}^2{\Sigma}_2
( (dx^6)^2 \rbrace,
$$
\be\label{9}
a_{ij} dx^i dx^j=\epsilon x^3g_{i_1 j_1}dx^{i_1}dx^{j_1}+2 g_{13} dx^2 dx^3+
4Ag_{22}(\epsilon x^1- A \Sigma_1)(dx^3)^2+
\ee
$$
(\tilde {\epsilon} x^6 +a)g_{i_2 j_2}dx^{i_2}dx^{j_2}+
2 g_{46} dx^5 dx^6+4\tilde{A}g_{55}(\tilde{\epsilon} x^4+\tilde{ A}
{\Sigma_1})
(dx^6)^2,
$$
\be\label{10}
\varphi=\frac{3}{2}(f_3+f_6)+c,
\ee
\be\label{11}
h_{ij}=a_{ij}+3(f_3+f_6+c)g_{ij},
\ee
where  $f_3=\epsilon x^3$, $f_6=\tilde{\epsilon}x^6+a$,
$\epsilon, \tilde{\epsilon}=0, 1$, $c$ and $a$ are some constants,
$a \ne 0$ when $\tilde{\epsilon}=0$,
$i_1, j_1=1,2, 3$, $i_2, j_2=4, 5, 6$,  $e_i =\pm 1$,
\be\label{12}
A=\epsilon x^2+\theta(x^3),
\quad
\tilde{A}=\tilde{\epsilon} x^4+\omega(x^6),
\quad
\Sigma_1=3 (f_6-f_3)^{-1},
\quad
\Sigma_2=3 (f_6-f_3)^{-2},
\ee
$\theta(x^3)$, $\omega(x^6)$ are arbitrary functions,
$\theta(x^3)\ne 0$ when $\epsilon=0$ and  $\omega(x^6) \ne 0$
when $\tilde \epsilon=0$.}

\section{Projective group characteristics of the
\\
6-dimensional type $[33]$ $h$-space}

To move further, we need a necessary and sufficient
condition of the constant curvature of the $h$-space under consideration.
This condition is known to be expressed by the formula
\be\label{13}
R_{jkl}^i=K({\delta_k}^i g_{jl}-{\delta_l}^i g_{jk}),
\quad
K={\rm const}.
\ee
Calculating the components of curvature tensor of our $h$-space,
we observe that the nonzero components
$R_{\nu b_p \mu}^{a_p}$
($b_p \ne \mu$, $a_p < b_p$) of the curvature tensor,
$R_{5 a_1 6}^{a_1}$, $R_{6 a_1 5}^{a_1}$,
$R_{6 a_1 6}^{a_1}$, $R_{2 a_2 3}^{a_2}$,
$R_{3 a_2 2}^{a_2}$, $R_{3 a_2 3}^{a_2}$  ($p=1, 2$, $a_1, b_1=1, 2, 3$,
$a_2, b_2=4, 5, 6$, $\nu$ and $\mu$ are the indices of nonzero
components of metric $g_{\nu \mu}$, (\ref{8})) are proportional to
$\epsilon$ or $\tilde \epsilon$.
In particulary, $R_{123}^2=\frac{3 {\epsilon}^2}{8 A}$,
$R_{456}^5=\frac{3 {\tilde\epsilon}^2}{8 \tilde  A}$.
Now one can derive from (\ref{13}) that $R_{123}^2=R_{163}^6$ and
$R_{456}^5=R_{436}^3$ which implies $\epsilon=\tilde\epsilon=0$.
Obviously, when $\epsilon=\tilde\epsilon=0$, all the
components of curvature tensor of the $h$-space we consider are equal
to zero.

Let us now prove the following theorems

\bigskip

\noindent
{\bf Theorem 2.} {\it The defining function of
projective motion in the type ${\rm [33]}$ $h$-space
of nonconstant curvature is $\phi=a_1 \varphi$,
where $a_1$ is a constant and
$\varphi$ is determined by {\rm (\ref{10})}.}

\noindent
{\bf Proof.} Suppose we are given with a vector field ${\xi}^i$
that gives the projective
motion with the defining function $\phi$ in the
$h$-space of the type $[33]$.
Then, for the tensor $b_{ij}={\nl{\xi}L} g_{ij}$ the following
Eisenhart equations are fulfilled
\be\label{14}
b_{ij,k}=2g_{ij} \phi_{,k}+g_{ik} \phi_{,j}+
g_{jk} \phi_{,i}.
\ee
Their integrability conditions are
\be\label{15}
b_{mi} R_{jkl}^m+b_{mj} R_{ikl}^m=g_{ik} \phi_{,jl}+g_{jk}\phi_{,il}-g_{li}
\phi_{,jk}-g_{lj} \phi_{,ik}.
\ee
Let us show that for any projective motion in the
type $[33]$ $h$-space the following conditions are satisfied
\be\label{16}
b_{\alpha \beta}=0,
\quad
\phi_{,\alpha \beta}=0,
\ee
$\alpha$, $\beta$ being indices of the components $g_{\alpha \beta}$
of the metric (\ref{8}) equal to zero. To prove that
$b_{11}=\phi_{,11}=0$, set $(ijkl)=(1212)$, $(1416)$ and $(1516)$
in (\ref{15}). Then, one gets
\be\label{17}
b_{11} R_{212}^1=-g_{22} \phi_{,11},
\quad
b_{11} R_{416}^1=-g_{46} \phi_{,11},
\ee
\be\label{18}
b_{11}(\frac{R_{212}^1}{g_{22}}-\frac{R_{516}^1}{g_{56}})=0.
\ee
Since $R_{416}^1=0$ and $g_{46} \ne 0$, we get $\phi_{,11}=0$.
Then, from (\ref{17}) it follows that $b_{11}=0$, since
$R_{212}^1 \ne 0$.
Similarly one can get the other equalities (\ref{16}).

Now, from (\ref{14}) with $(ijk)=(166),$ (266), (433), (533), (346) and
(613) and taking into account (\ref{16}), one obtains
$$
\phi_{,1}=\phi_{,2}=\phi_{,4}=\phi_{,5}=0,
$$
$$
\phi_{,3}=\frac{3}{2}\frac{\epsilon}{f_6-f_3}
(\frac{b_{46}}{g_{46}}-\frac{b_{13}}{g_{13}}),
\quad
\phi_{,6}=\frac{3}{2}\frac{\tilde\epsilon}{f_6-f_3}
(\frac{b_{46}}{g_{46}}-\frac{b_{13}}{g_{13}}).
$$
Comparing the last two formulas one finds
\be\label{19}
\tilde\epsilon \phi_{,3}=\epsilon \phi_{,6}.
\ee
Then, from the equality $\phi_{,36}=0$ it follows that $\pl_{36} \phi=0$.
If $\tilde\epsilon \ne 0$, differentiating (\ref{19}) w.r.t. $x^3$
gives  $\pl_{33} \phi=0$. If $\tilde\epsilon=0$, then
$\epsilon \ne 0$, since we consider the space of nonzero curvature.
In the case of $\tilde\epsilon=0$ and $\epsilon \ne 0$, the relation
$\pl_{33}\phi=0$ is obtained from (\ref{15}) with
$(ijkl)=(3312)$, (3436). One also can analogously show that
$\pl_{66} \phi=0$. Integrating the equation obtained, one finally finds
\be\label{20}
\phi=\frac{1}{2} a_1 \nl{i=1}{\sum}^6 f_i=a_1 \varphi,
\quad
f_1=f_2=f_3=\epsilon x^3,
\quad
f_3=f_4=f_6=\tilde\epsilon x^6+a,
\ee
where $\varphi$ is determined by formula (\ref{10}).

\bigskip

\noindent
{\bf Theorem 3.} {\it Any covariant constant symmetric
tensor $b_{ij}$ in the 6-dimensional
type $[33]$ $h$-space of nonzero curvature is proportional to the
fundamental tensor
\be\label{21}
b_{ij}=a_2 g_{ij},
\quad
(a_2={\rm const}).
\ee }

\noindent
{\bf Proof.} Let $b_{ij}$ in $V^6$ satisfy the equation
\be\label{22}
b_{ij,k}=0.
\ee
Integrability conditions for this equation follow from (\ref{15})
with $\phi=const$. Then, from (\ref{16}) one has
\be\label{23}
b_{\alpha \beta}=0.
\ee
Note that for the proof of the theorem one suffices  to show that
$b_{ij}=a_2 g_{ij}$, where $a_2$ is a coefficient which is constant,
due to the equation $b_{ij,k}=0$ and parallelism of the metric tensor.

Putting in equation (\ref{15}) $(ijkl)=(2223)$, (3312), (3331),
one finds
\be\label{24}
b_{22} R_{223}^2-b_{23} R_{232}^3=0,
\quad
b_{13} R_{312}^1-b_{23} R_{232}^3=0,
\ee
$$
b_{13} R_{331}^1+b_{23} R_{331}^2+b_{33} R_{331}^3=0.
$$
From (\ref{22}) with  $(ijk)=(436)$, (535). (536), (636)
it follows that
\be\label{25}
f'_{3}(\frac{b_{33}}{g_{33}}-\frac{b_{46}}{g_{46}})=0,
\quad
f'_{3}(\frac{b_{33}}{g_{33}}-\frac{b_{55}}{g_{55}})=0,
\\
f'_{3}(\frac{b_{33}}{g_{33}}-\frac{b_{56}}{g_{56}})=0,
\quad
f'_{3}(\frac{b_{33}}{g_{33}}-\frac{b_{66}}{g_{66}})=0.
\ee
If $f'_{3} \ne 0$, one obtains (\ref{21})
from (\ref{24}), (\ref{25}) and (\ref{23}).
Otherwise, the wanted result follows from
(\ref{22}) with $(ijk)=(163)$, (262), (263), (363) and from
(\ref{15}) with $(ijkl)=(5556)$, (6645), (6664).

From this theorem, one immediately obtains

\bigskip

\noindent
{\bf Theorem 4.} {\it Affine group of the 6-dimensional
type $[33]$ $h$-space of non-constant curvature consists of
the homothetics. }

This theorems and linearity of the Eisenhart
equation give the general solution of the Eisenhart equation
in the type $[33]$ $h$-space of non-constant curvature
in the form  $a_1 h_{ij}+a_2 g_{ij}$ with two arbitrary constants
$a_1, a_2$,
$$
h_{ij}=a_{ij}+(\sum_{k=1}^6 f_k) g_{ij},
$$
tensors $g_{ij}$, $a_{ij}$ being determined by (\ref{8})
and (\ref{9}).

Hence, one obtains

\bigskip

\noindent
{\bf Theorem 5.} {\it All projective motion of the
{\rm 6}-dimensional type $[33]$ $h$-space of non-constant curvature
are obtained by integrating the equation
$$
L_{\xi} g_{ij}=\xi_{i,j}+\xi_{j,i}=a_1 a_{ij}+(a_1 \sum_{k=1}^6 f_k+a_2)g_{ij},
$$
where tensors $g_{ij}$ and  $a_{ij}$ are determined by
{\rm (\ref{8})} and {\rm (\ref{9})}, and $a_1$, $a_2$ are constants.}

It leads to the following important group characteristic of the
type $[33]$ $h$-spaces:

\bigskip

\noindent
{\bf Theorem 6.} {\it If the type $[33]$ $h$-space
of non-constant curvature admits a nonhomothetic projective Lie algebra $P_r$,
then this algebra contains the subalgebra $H_{r-1}$ of infinitesimal
homothetics of dimension $r-1$.}

\bigskip

I am grateful to A.V.Aminova for the constant encouragement and
discussions. I also acknowledge Koji Matsumoto for the discussions and kind
hospitality at the Yamagata University where this paper was completed.
The work is partially supported by RFBR grant 01-02-17682a and INTAS grant 00-00334.

\end{document}